\documentclass[12pt]{article}
\usepackage{Wsom97}
\begin{document}
\title{\bf Theoretical Aspects of the SOM Algorithm}
%
\author{\large M.Cottrell$^\dagger$,
J.C.Fort$^\ddagger$, G.Pag\`es$^\ast$\\
  \normalsize $\dagger$ SAMOS/Universit\'e Paris 1\\
  \normalsize 90, rue de Tolbiac, F-75634 Paris Cedex 13, France\\
  \normalsize Tel/Fax : 33-1-40-77-19-22, E-mail: cottrell@univ-paris1.fr\\
  \normalsize $\ddagger$ Institut Elie Cartan/Universit\'e Nancy 1 et SAMOS\\
  \normalsize F-54506 Vand\oe uvre-L\`es-Nancy Cedex, France\\
  \normalsize E-mail: fortjc@iecn.u-nancy.fr\\
  \normalsize $\ast$ Universit\'e Paris 12 et Laboratoire de Probabilit\'es
                        /Paris 6\\
  \normalsize F-75252 Paris Cedex 05, France\\
  \normalsize E-mail:gpa@ccr.jussieu.fr}

\date{}

\newtheorem{madef}{Definition}
\newtheorem{monthe}{Theorem}

\maketitle
\thispagestyle{empty}
\begin{abstract}
The SOM algorithm is very astonishing. On the one hand, it is very simple to 
write down and to simulate, its practical properties are clear and easy 
to observe. But, on the other hand, its theoretical properties still 
remain without proof in the general case, despite the great efforts
of several authors. In this paper, we pass in review the last results 
and provide some conjectures for the future work. 

\end{abstract}

\section{Introduction}

The now very popular SOM algorithm was originally devised by Teuvo 
Kohonen in 1982 \cite{KOHO1} and \cite{KOHO2}. It was presented as a model 
of the self-organization of neural connections. What immediatly raised the 
interest of the scientific community (neurophysiologists, computer 
scientists, mathematicians, physicists) was the ability of such a simple 
algorithm to produce organization, starting from possibly total disorder. 
That is called the {\em self-organization} property.

As a matter of fact, the algorithm can be considered as a generalization of 
the Competitive Learning, that is a Vectorial Quantization 
Algorithm \cite{IEEE}, without any notion of neighborhood between the units.  

In the SOM algorithm, a neighborhood structure is defined for the units 
and is respected  throughout the learning process, which imposes the 
conservation of the neighborhood relations. So the weights are 
progressively updated according to the presentation of the inputs, in such 
a way that neighboring inputs are little by little mapped onto the same 
unit or neighboring units. 

There are two phases. As well in the practical applications as in the 
theoretical studies, one can observe self-organization first (with large 
neighborhood and large adaptation parameter), and later on convergence 
of the weights in order to quantify the input space. In this second phase, 
the adaptation parameter is decreased to 0, and the neighborhood is small 
or indeed reduced to one unit (the organization is supposed not to be 
deleted by the process in this phase, that is really true for the 0-neighbor
setting).

Even if the properties of the SOM algorithm can be easily reproduced 
by simulations, and despite all the efforts, the Kohonen algorithm is 
surprisingly resistant to a complete mathematical study. As far as we know, 
the only case where a complete analysis has been achieved is the 
{\em one dimensional case} (the input space has dimension 1) for a 
{\em linear network} (the units are disposed along a one-dimensional array). 

A sketch of the proof was provided in the Kohonen's original papers 
\cite{KOHO1}, \cite{KOHO2} in 1982 and in his books \cite{KOHO3}, 
\cite{KOHO6} in 1984 and 1995. The  first complete proof of both 
self-organization and convergence properties was established (for uniform 
distribution of the inputs and a simple step-neighborhood function) by 
Cottrell and Fort in 1987, \cite{COTT}. 

Then, these results were generalized to a wide class of input distributions  
by Bouton and Pag\`es in 1993 and 1994, \cite{BOUT2}, \cite{BOUT3} and to a 
more general neighborhood by Erwin et al. (1992) who have sketched the 
extension of the proof of self-organization \cite{ERWI2} and studied the 
role of the neighborhood function \cite{ERWI1}. Recently, Sadeghi 
\cite{SADE1}, \cite{SADE2} has studied the self-organization for a general type of stimuli distribution and neighborhood function.

At last, Fort and Pag\`es in 1993, \cite{FORT4}, 1995 \cite{FORT5}, 
1997 \cite{FORT10}, \cite{FORT11} (with Benaim) have achieved the 
rigorous proof of the 
almost sure convergence towards a unique state, after self-organization, 
for a very general class of neighborhood functions.

Before that, Ritter et al. in 1986 and 1988, \cite{RITT1}, \cite{RITT2} 
have thrown some light on the stationary state in any dimension, but they 
study only the final phase {\em after the self-organization}, and do not 
prove the existence of this stationary state.

\medskip

In multidimensional settings, it is not possible to define what could be 
a {\em well ordered configuration set} that would be stable for the 
algorithm and that could be an absorbing class. For example, the 
grid configurations that Lo et al. proposed in 1991 or 1993, \cite{LO1}, 
\cite{LO2} are  not stable as proved in \cite{COTT1}. Fort and Pag\`es in 
1996, \cite{FORT6} show that 
there is no organized absorbing set, at least when the stimuli space is 
continuous. On the other hand, Erwin et al. in 1992 
\cite{ERWI2} have proved that it is {\em impossible to associate a global 
decreasing potential function} to the algorithm, as long as the probability 
distribution of the inputs is continuous. Recently,  Fort and Pag\`es 
in 1994, \cite{FORT4}, in 1996 \cite{FORT5} and \cite{FORT6}, 
Flanagan in 1994 and 1996 \cite{FLAN}, \cite{FLAN1} 
gave some results in higher dimension, but these remain incomplete. 

\medskip

In this paper, we try to present the state of the art. As a continuation 
of previous paper \cite{COTT5}, we gather the more recent 
results that have been published in different journals that can be not 
easily get-a-able for the neural community.
 
We do not speak about the variants of the algorithm 
that have been defined and studied by many authors, in order to improve 
the performances or to facilitate the mathematical analysis, see for 
example \cite{BISH}, \cite{LUTT}, \cite{RUZI}, \cite{THI1}. We do not 
either address the numerous applications of the SOM 
algorithm. See for example the Kohonen's book \cite{KOHO6} to have an 
idea of the profusion of these applications. We will only mention as a 
conclusion some original data analysis methods based on the SOM algorithm.


The paper is organized as follows: in section 2, we define the notations.
The section 3 is devoted to the one dimensional case. Section 4 deals 
with the multidimensional 0-neighbor case, that is the simple competitive 
learning and gives some light on the quantization performances. 
In section 5, some partial results about the multidimensional setting are 
provided. Section 6 treats the discrete finite case and we present 
some data analysis methods derived from the SOM algorithm. The conclusion 
gives some hints about future researches.

\section{Notations and definitions}

The network includes $n$ units located in an ordered lattice (generally in 
a one- or two-dimensional array). If $I=\{1, 2, \ldots, n\}$ is the set 
of the indices, the neighborhood structure is provided by a neighborhood 
function $\Lambda$ defined on $I \times I$. It is symmetrical,
non increasing, and depends only on the distance between $i$ and 
$j$ in the set of units $I$, (e.g. $\mid i-j \mid$ if 
$I=\{1, 2, \ldots, n\}$ is one-dimensional). $\Lambda (i,j)$ decreases 
with increasing distance between $i$ and $j$, and $\Lambda (i,i)$ is usually equal to 1. 

\medskip

The input space $\Omega$ is a bounded convex subset of ${\cal{R}}^{d}$, 
endowed with the Euclidean distance. The inputs $x(t), t\geq 1$ are 
$\Omega$-valued, independent with common distribution $\mu$.

The network state at time $t$ is given by 
\[m(t) = (m_1(t),\:m_2(t),\:\ldots,\:m_n(t)).\]
where $m_i(t)$ is the $d$-dimensional weight vector of the unit $i$. 

For a given state $m$ and input $x$, the {\bf winning} unit 
$i_c(x, m)$ is the unit whose weight $m_{i_c(x, m)}$ 
is the closest to the input $x$. Thus the network defines a map 
$\Phi_{m}: x \longmapsto i_c(x, m)$, from $\Omega$ to $I$, and 
the goal of the learning algorithm is to converge to a network state such  
the $\Phi_{m}$ map will be ``topology preserving''in some sense.

For a given state $m$, let us denote $C_i(m)$ the set of the inputs such 
that $i$ is the winning unit, that is $C_i(m)=\Phi_{m}^{-1}(i)$. 
The set of the classes $C_i(m)$ is the Euclidean Vorono\"{\i} tessellation 
of the space $\Omega$ related to $m$.

\medskip

The SOM algorithm is recursively defined by :

\begin{equation}
\left\{
\begin{array}{lll}
i_c(x(t+1),m(t)) & = &
            \mbox{argmin }\left\{ \|x(t+1)-m_i(t)\|,i \in I \right\} \\
m_i(t+1) & = & m_i(t)-\varepsilon_t \Lambda(i_{0},i) (m_i(t)-x(t+1)), 
                                                    \forall i \in I 
\end{array}
\right.
\end{equation}

\medskip
\medskip

The essential parameters are 
\begin{itemize}
\item the dimension $d$ of the input space
\item the topology of the network
\item the adaptation gain parameter $\varepsilon_t$, 
which is $]0,1[$-valued, constant or decreasing with time, 
\item the neighborhood function $\Lambda$, which can be constant or time 
dependent, 
\item the probability distribution $\mu$.
\end{itemize}

\medskip

{\bf Mathematical available techniques}

\medskip

As mentioned before, when dealing with the SOM algorithm, one has to 
separate two kinds of results: those related to self-organization, and 
those related to convergence after organization. In any case, all the 
results have been obtained for a fixed time-invariant neighborhood function.

First, the network state at time $t$ is a random 
$\Omega^n$-valued vector $m(t)$ displaying as~:

\begin{equation}
m(t+1) = m(t) - \varepsilon_t \: H(x(t+1),\,m(t))  \label{defi}
\end{equation}
(where $H$ is defined in an obvious way according to the updating equation)
is a stochastic process. If $\varepsilon_t$ and $\Lambda$ are time-invariant,
it is an homogeneous {\em Markov chain} and can be studied with the usual 
tools if possible (and fruitful). For example, if the algorithm converges 
in distribution, this limit distribution has to be an invariant measure 
for the Markov chain. If the algorithm has some fixed point, this point 
has to be an absorbing state of the chain. If it is possible to prove some 
strong organization \cite{FORT6}, it has to be associated to an absorbing 
class.

Another way to investigate self-organization and convergence is to study the 
associated ODE (Ordinary Differential Equation) \cite{KUSH} that describes the 
mean behaviour of the algorithm :

\begin{equation}
\frac{dm}{dt} \, = \, - \:\: h(m)
\label{ODE}
\end{equation}
 
where

\begin{equation}
h(m)=E(H(x,\,m)) = \int H(x,\,m) \: d\mu(x)
\end{equation}

is the expectation of $H(.,\,m)$ with respect to the probability measure $\mu$. 

Then it is clear that all the possible limit states $m^{\star}$ are solutions of the functional equation 
$$h(m)=0$$

and any knowledge about the possible attracting equilibrium points of the ODE 
can give some light about the self-organizing property and the 
convergence. But actually the complete asymptotic study of the ODE in the 
multidimensional setting seems to be untractable. One has to verify some 
{\em global assumptions} on the function $h$ (and on its {\em gradient}) 
and the explicit calculations are quite difficult, and perhaps impossible. 

\medskip

In the {\em convergence phase}, the techniques depend on the kind of the 
desired convergence mode. For the {\em almost sure} convergence, the 
parameter $\varepsilon_t$ needs to decrease to 0, and the form of 
equation (\ref{defi}) suggests to consider the SOM algorithm as a 
Robbins-Monro \cite{ROBB} algorithm.

The usual hypothesis on the adaptation parameter to get almost sure 
results is then:

\begin{equation}
        \sum_t \varepsilon_t = +\infty \: \mbox{ and } 
        \sum_t \varepsilon_t^2 < + \infty.                 
   \label{epsi}
\end{equation}

The less restrictive conditions $\sum_t \varepsilon_t = +\infty$ and 
$\varepsilon_t \searrow 0$ generally do not ensure the almost sure 
convergence, but some weaker convergence, for instance the convergence in 
probability. 

Let us first examine the results in dimension 1.

\section{The dimension 1}

\subsection{The self-organization}

The input space is $[0,1]$, the dimension $d$ is 1 and
the units are arranged on a linear array. The neighborhood function $\Lambda$ 
is supposed to be {\em non increasing} as a function of the distance between 
units, the classical step neighborhood function satisfies this condition. 
The input distribution $\mu$ is {\em continuous} on $[0,1]$: this means 
that it does not weight any point. This is satisfied for example by any 
distribution having a {\em density}.

Let us define 
$$F_n^+ = \{m \in {\cal {R}} \, / \, 0 < m_1 < m_2 < \ldots < m_n <1 \}$$ 
and
$$F_n^- = \{m \in {\cal {R}} \, / \, 0 < m_n < m_{n-1} <\ldots < m_1<1 \}.$$

In \cite{COTT}, \cite{BOUT2}, the following results are proved using 
Markovian methods :

\begin{monthe}

(i) The two sets $F_n^+$ and $F_n^-$ are absorbing sets.

\noindent (ii) If $\varepsilon$ is constant, and if $\Lambda$ is decreasing 
as a function of the distance (e.g. if there are only two neigbors) the 
entering time $\tau$,  that is the hitting time of $F_n^+ \cup F_n^-$, is 
almost surely finite, and $\exists \lambda > 0, \mbox{ s.t. } 
\sup_{ m \in [0,1]^n} \: E_{m} (\exp ( \lambda \tau))$ is finite, 
where $E_{m}$ denote the expectation given ${\bf m}(0)=m$.
\label{orga1}

\end{monthe}

The theorem \ref{orga1} ensures that the algorithm will almost surely 
order the weights. These results can be found for the more particular case 
($\mu$ uniform and two neighbors) in Cottrell and Fort \cite{COTT}, 1987, 
and the succesive generalisations in Erwin et al. \cite{ERWI2}, 1992, 
Bouton and Pag\`es \cite{BOUT2}, 1993, Fort and Pag\`es \cite{FORT5}, 1995, 
Flanagan \cite{FLAN1}, 1996.

The techniques are the Markov chain tools.

Actually following \cite{BOUT2}, it is possible to prove that 
whenever $\varepsilon \searrow 0$ and $\sum \varepsilon_t = + \infty$, then 
$\forall m \in [0,1]^n, \mbox{Proba}_m ( \tau < + \infty ) > 0$, (that is the probability of self-organization is positive regardless the initial values, but not a priori equal to 1).  In \cite{SADE2}, Sadeghi uses a generalized definition of the winner unit and shows that the probability of 
self-organization is uniformly positive, without assuming a lower bound for 
$\varepsilon_t$.

No result of almost sure reordering with a vanishing $\varepsilon_t$ is known 
so far. In \cite{COTT1}, Cottrell and Fort propose a still not proved 
conjecture: it seems that the re-organization occurs when the parameter 
$\varepsilon_t$ has a $\frac{1}{\ln t}$ order.

\subsection{The convergence for dimension 1}

After having proved that the process enters an ordered state set (increasing 
or decreasing), with probability 1, it is possible to study the convergence 
of the process. So we assume that ${\bf m}(0) \in F_n^+$. It would be the same if 
${\bf m}(0) \in F_n^-$.

\subsubsection{Decreasing adaptation parameter}

In \cite{COTT} (for the uniform distribution), in \cite{BOUT3}, \cite{FORT5} 
and more recently in \cite{FORT10}, \cite{FORT11}, 1997, the almost sure 
convergence is proved in a very general setting. The results are gathered 
in the theorem below :

\begin{monthe}

Assume that \\
\noindent 1) $(\varepsilon_t) \in ]0, 1[$ satisfies the condition 
(\ref{epsi}), 

\noindent 2) the neighborhood function satisfies the condition 
$H_{\Lambda}$: there exists $k_0 < \frac{n-1}{2}$ such that 
$\Lambda(k_0 +1) < \Lambda(k_0)$,


\noindent 3) the input distribution $\mu$ satisfy the condition $H_\mu$: 
it has a density $f$ such that $f>0$ on $]0,1[$  and  $\ln (f)$ is strictly 
concave (or only concave, with  $\lim_{0^+} f  + \lim_{1^-} f$ positive),

\noindent Then

\noindent (i) The mean function $h$ has a unique zero $m^{\star}$ in 
$F_n^+$.

\noindent (ii) The dynamical system $\frac{dm}{dt} = - h(m)$ is 
cooperative on $F_n^+$, i.e. the non diagonal elements of 
$\nabla h(m)$ are non positive.

\noindent (iii) $m^{\star}$ is attracting.

So if ${\bf m}(0) \in F_n^+$, ${\bf m}(t) \stackrel{a.s} {\longrightarrow}~m^{\star}$ 
almost surely.

\label{cvdim1}

\end{monthe}

In this part, the authors use the ODE method, a result by M.Hirsch on 
cooperative dynamical system \cite{HIRS}, and the Kushner \& Clark 
Theorem \cite{KUSH}, \cite{FORT10}. A.Sadeghi put in light that the 
non-positivity of non-diagonal terms of $\nabla h$ is exactly the basic 
definition of a cooperative dynamical system and he obtained partial 
results in \cite{SADE1} and more general ones in \cite{SADE2}.

\medskip

We can see that the assumptions are very general. Most of the usual 
probability distributions (truncated on $[0,1]$) have
a density $f$ such that  $\ln (f)$  is strictly concave. 
On the other hand, the uniform distribution is not strictly ln-concave as 
well as the truncated exponential distribution, but both cumply the 
condition $\lim_{0^+} f  + \lim_{1^-} f$ positive. 

Condition (\ref{epsi}) is essential, because if  $\varepsilon_t  
\searrow 0$ and $\sum_t \varepsilon_t = + \infty$, there is only a priori 
convergence in probability.

In fact, by studying the associated ODE, Flanagan \cite{FLAN} 
shows that before ordering, it can appear metastable equilibria.

In the uniform case, it is possible to calculate the limit $m^{\star}$. 
Its coordinates are solutions of a $(n\times n)$-linear system which 
can be found in \cite{KOHO3} or \cite{COTT}. An explicit expression,
up to the solution of a $3 \times 3$ linear system is proposed in 
\cite{BOUT2}. Some further investigations are made in \cite{FORT9}.

\subsubsection{Constant adaptation parameter}

Another point of view is to study the convergence of ${\bf m}(t)$ when 
$\varepsilon_t = \varepsilon$ is a constant. Some results are available 
when the neighborhood function corresponds to the two-neighbors setting. 
See \cite{COTT}, 1987, (for the uniform distribution) and \cite{BOUT3}, 
1994, for the more general case. One part of the results also hold for a 
more general neighborhood function, see \cite{FORT10}, \cite{FORT11}.

\begin{monthe} Assume that ${\bf m}(0) \in F^+_n$,

\noindent Part A: Assume that the hypotheses $H_{\mu}$ and $H_{\Lambda}$ 
hold as in Theorem \ref{cvdim1}, then

\noindent For each $\varepsilon \in ]0,1[$, there exists some
invariant probability $\nu^{\varepsilon}$ on $F_n^+$.

\noindent Part B: Assume only that $\Lambda(i,j)=1$ if and only if $|i-j|=0$ 
or $1$ (classical 2-neighbors setting),
 
\noindent(i) If the input distribution $\mu$ has an absolutely continuous 
part (e.g. has a density), then for each $\varepsilon \in ]0, 1[$, there 
exists a unique probability distribution $\nu^{\varepsilon}$ such that 
the distribution of ${\bf m}^t$ weakly converges to $\nu^{\varepsilon}$ when 
$t \longrightarrow \infty$. The rate of convergence is geometric. 
Actually the Markov chain is Doeblin recurrent.

\noindent (ii) Furthermore, if $\mu$ has a positive density, $\forall 
\varepsilon$, $\nu^{\varepsilon}$ is equivalent to the Lebesgue measure 
on $F^+_n$ if and only if $n$ is congruent with 0 or 1 modulo 3. If $n$ is 
congruent with 2 modulo 3, the Lebesgue measure is absolutely continuous 
with respect to $\nu^{\varepsilon}$ , but the inverse is not true, that is 
$\nu^{\varepsilon}$ has a singular part.  

\noindent Part C: With the general hypotheses of Part A (which includes 
that of Part B), if $m^{\star}$ is the unique globally attractive 
equilibrium of the ODE (see Theorem \ref{cvdim1}), thus $\nu^{\varepsilon}$ 
converges to the Dirac distribution on $m^\star$ when 
$\varepsilon \searrow 0$ .
\label{eps_constant}
\end{monthe}

So when $\varepsilon$ is very small, the values will remain very close 
to $m^{\star}$.

Moreover, from this result we may conjecture that for a suitable choice of 
$\varepsilon_t$, certainly $\varepsilon_t = \frac{A}{\ln t}$, where $A$ is a constant, both 
self-organization and convergence towards the unique $m^{\star}$ can be 
achieved. This could be proved by techniques very similar to the 
simulated annealing methods.

\section{The 0 neighbor case in a multidimensional setting}

In this case, we take any dimension $d$, the input space is $\Omega \subset 
{\cal R}^d$ and $\Lambda(i,j)=1$ if $i=j$, and $0$ elsewhere. There is no 
more topology on $I$, and {\em reordering} no makes sense. In this case the 
algorithm is essentially a stochastic version of the Linde, Gray and Buzo 
\cite{LIND} algorithm (LBG). It belongs to the family of the vectorial 
quantization algorithms and is equivalent to the Competitive Learning. The 
mathematical results are more or less reachable. Even if this algorithm is 
deeply different from the usual Kohonen algorithm, it is however interesting 
to study it because it can be viewed as a limit situation when the 
neighborhood size decreases to 0. 

\medskip

The first result (which is classical for Competitive learning), and can be 
found in \cite{RITT3}, \cite{PAGE1}, \cite{KOHO5} is: 

\begin{monthe}

(i) The 0-neighbor algorithm derives from the potential

\begin{equation}
V_n(m)=\frac{1}{2}\int \min_{1\leq i\leq n}\| m_i-x \|^2 d{\mu}(x)
\label{vari1}
\end{equation}

\noindent (ii) If the distribution probability $\mu$ is continuous (for 
example $\mu$ has a density $f$), 
\begin{equation}
V_n(m)=\frac{1}{2}\sum_{i=1}^{n}\int_{C_{i}(m)} 
\|m_i-x\|^2 f(x)dx
= \frac{1}{2}\int \min_{1\leq i \leq n} \|m_i-x\|^2 f(x)dx 
\label{vari2}
\end{equation}
where $C_i(m)$ is the Vorono\"{\i} set related with the unit $i$ for the 
current state $m$. 
\label{potential}
\end{monthe}

The potential function $V_n(m)$ is nothing else than the {\em intra-classes 
variance} used by the statisticians to characterize the quality of a 
clustering. In the vectorial quantization setting, $V_n(m)$ is called 
{\em distortion}. It is a measure of the loss of information when replacing 
each input by the closest weight vector (or {\em code vector}). The potential 
$V_n(m)$ has been extensively studied since 50 years, as it can be seen in 
the Special Issue of IEEE Transactions on Information Theory (1982), 
\cite{IEEE}.

\medskip

The expression (\ref{vari2}) holds as soon as $m_i \neq m_j$ for all 
$i \neq j$ 
and as the borders of the Vorono\"{\i} classes have probability 0, 
($\mu (\cup_{i=1}^n \partial C_i(m))=0$). This last condition is always 
verified when the distribution $\mu$ has a density $f$. With these two 
conditions, $V(m)$ is  differentiable at $m$ and its gradient vector reads

$$\nabla V_n(m) = \left(\int_{C_i(m)}(m_i - x) f(m) d(m)\right).$$

So it becomes clear (\cite{PAGE1},\cite{KOHO6}) that the Kohonen algorithm 
with 0 neighbor is the stochastic gradient descent relative to the function 
$V_n(m)$ and can be written :

$$m(t+1) = m(t) - \varepsilon_{t+1} 
{\bf 1}_{C_i(m(t))}(x(t+1))(m(t) - x(t+1))$$

\noindent where ${\bf 1}_{C_i(m(t))}(x(t+1))$ is equal to 1 if 
$x(t+1) \in C_i(m(t))$, and 0 if not.

\medskip

The available results are more or less classical, and  can be found in 
\cite{LIND} and \cite{BOUT4}, for a general dimension $d$ and 
a distribution $\mu$ satisfying the previous conditions. 

Concerning the convergence results, we have the following when the dimension 
$d=1$, see Pag\`es (\cite{PAGE1}, \cite{PAGE2}), the Special Issue in 
IEEE \cite{IEEE} and also \cite{LAMB} for (ii):

The parameter $\varepsilon(t)$ has to satisfy the conditions (\ref{epsi}).

\begin {monthe} {\bf Quantization in dimension 1}\\
(i) If $\nabla V_n$ has finitely many zeros in $F_n^+$, $m(t)$ converges 
almost surely to one of these local minima.

\noindent (ii) If the hypothesis $H_{\mu}$ holds (see Theorem (\ref{cvdim1})), 
$V_n$ has only one zero point in $F^+_n$, say $m^{\star}_n$. This point 
$m^{\star}_n \in F^+_n$ and is a minimum. Furthermore if $m(0) \in F_n^+$, 
$m(t) \stackrel {a.s.} {\longrightarrow} m^{\star}_n$.

\noindent (iii) If the stimuli are uniformly distributed on $[0,1]$, then 
$$ m^{\star}_n\,=\, ((2i-1)/2n)_{1 \leq i \leq n}.$$

\label{Quantdim1}

\end{monthe}

The part (ii) shows that the global minimum de $V_n(m)$ is reachable in the 
one-dimensional case and the part {\em (iii)} is a confirmation of the fact 
that the algorithm provides an optimal discretization of continous 
distributions. 

\medskip
A weaker result holds in the $d$-dimensional case, because one has only 
the convergence to a {\em local minimum} of $V_n(m)$. 
\medskip
 
\begin {monthe} {\bf Quantization in dimension d}\\
If $\nabla V_n$ has finitely many zeros in $F_n^+$, and if these zeros 
have all their components pairwise distinct, $m(t)$ converges almost surely 
to one of these local minima.

\label{Quantdimd}
\end{monthe}

In the $d$-dimensional case, we are not able to compute the limit, even in 
the uniform case. Following \cite{NEWN} and many experimental results, 
it seems that the minimum distortion could be reached for an hexagonal 
tesselation, as mentioned in \cite{FORT9} or \cite{KOHO6}.

\medskip

In both cases, we can set the properties of the global minima of $V_n(m)$, in 
the general $d$-dimensional setting. Let us note first that $V_n(m)$ is 
invariant under any permutation of the integers $1,2,\ldots,n$. So we can 
consider one of the global minima, the ordered one  (for example the
lexicographically ordered one).

\begin{monthe}{\bf Quantization property}\\
(i) The function $V_n(m)$ is continuous on $({\cal R}^d)^n$ and 
reaches its (global) minima inside $\Omega^n$.

\noindent (ii) For a fixed $n$, a point $m^{\star}_n$ at which the function 
$V_n$ is minimum has pairwise distinct components.

\noindent (iii) Let $n$ be a variable and 
$m^{\star}_n=(m^{\star}_{n,1},m^{\star}_{n,2},\ldots,m^{\star}_{n,n})$  
the ordered minimum of $V_n(m)$. The sequence 
$\min _{({\cal R}^d)^n}V_n(m)=V_n(m^{\star}_n)$ 
converges to 0 as $n$ goes to $+\infty$. 

More precisely, there exists a speed $\beta=2/d$ and a constante $A(f)$ 
such that 
$$n^{\beta}V_n(m^{\star}_n) \longrightarrow A(f)$$
when $n$ goes to $+\infty$.

Following Zador \cite{ZADO}, the constant 
$A(f)$ can be computed, $A(f)=a_{d}\parallel f \parallel_{\rho}$, 
where $a_d$ does not depend on $f$, $\rho=d/(d+2)$ and 
$\parallel f \parallel_{\rho} = [\int f^{\rho}(x) dx ]^{1/\rho}$.

\noindent (iv) Then, the weighted empirical discrete probability measure
$$\mu_n~ =~\sum_{i=1}^{n} 
{\mu}(C_i(m^{\star}_{n})) \delta_{m^{\star}_{n,i}}$$
converges in distribution to the probability measure ${\mu}$, when 
$n \rightarrow \infty$. 

\noindent (v) If $F_n$ (resp. $F$) denotes the distribution function of 
$\mu_n$ (resp. $\mu$), one has 
$$\min_{({\cal R}^d)^n} V_n(m)=
\min_{({\cal R}^d)^n} \int_{\Omega}(F_n(x) - F(x))^2 dx,$$
so when $n \rightarrow \infty$, $F_n$ converges to $F$ in quadratic norm. 
\label{VQD}
\end{monthe}

\medskip

The convergence in {\em (iv)} properly defines the {\em quantization 
property}, and explains how to reconstruct the input distribution from 
the $n$ code vectors after convergence.  But in fact this convergence 
holds for any sequence $y^{\star}_n = y_{1,n},y_{2,n},\ldots,y_{n,n}$, 
which ``fills '' the space when $n$ goes to $+\infty$: for example it is 
sufficient that for any $n$, there exists an integer $n' > n$ such that 
in any interval $y_{i,n},y_{i+1,n}$ (in ${\cal R}^d$), there are some points 
of $y^{\star}_{n'}$. But for any sequence of quantizers satisfying this 
condition, even if there is convergence in distribution, even if the speed 
of the convergence can be the same, the constant $A(f)$ will differ since 
it will not realize the minimum of the distortion.

For each integer $n$, the solution $m^{\star}_n$ which minimizes the 
quadratic distortion $V_n(m)$ and the quadratic norm  
$\parallel F_n - F \parallel^2$ is said to be {\em an optimal 
$n$-quantizer }. It ensures also that the discrete distribution function 
associated to the minimum $m^{\star}_n$ suitably weighted by the probability 
of the Vorono\"{\i} classes, converges to the initial distribution function 
$F$. So the 0-neighbor algorithm provides a skeleton of the input 
distribution and as the distortion tends to 0 as well as the quadratic norm 
distance of $F_n$ and $F$, it provides an {\em optimal quantizer}. 
The weighting of the Dirac functions by the volume of the Vorono\"{\i} 
classes implies that the distribution $\mu_n$ is usually quite different 
from the empirical one, in which each term would have the same weight $1/n$.

This result has been used by Pag\`es in \cite{PAGE1} and \cite{PAGE2} to 
numerically compute integrals. He shows that the speed of convergence of 
the approximate integrals is exactly $n^{\frac{2}{d}}$ for smooth enough 
functions, which is faster than the Monte Carlo method while $d \leq 4$. 

The difficulty remains that the optimal quantizer $m^{\star}_n$ is not 
easily reachable, since the stochastic process $m(t)$ converges only to a 
local minimum of the distortion, when the dimension is greater than 1.
 
\medskip

{\bf Magnification factor}

There is {\em some confusion} \cite{KOHO3}, \cite{RITT1}, between the asymptotic 
distribution of an {\em optimal quantizer} $m^{\star}_n$ when 
$n \longrightarrow \infty$ and that one of the best {\em random quantizer}, 
as defined by Zador \cite{ZADO} in 1982.  

The Zador's result, extended to the multi-dimensional case, is as follows :
{\em Let $f$ be the input density of the measure $\mu$, and $(Y_1, Y_2, 
\ldots, Y_n)$ a {\em random quantizer}, where the code vectors $Y_i$ are 
independent with common distribution of density $g$. 

Then, with some weak assumptions about $f$ and $g$, the distortion tends to 0 
when $n \longrightarrow \infty$, with  speed $\beta=2/d$, and it is possible 
to define the quantity} $$A(f,g) = \lim_{n \longrightarrow \infty} n^{\beta} 
     E_g [ \sum _{i=1}^{n} \int_{C_i} \|Y_i - x\|^2 
					f(x) dx ]$$

{\em Then for any given input density $f$, the density $g$ (assuming some 
weak condition) which minimises $A(f,g)$ is} 
$$g^{\star}~\sim C~f^{d/d+2}.$$

\medskip

The inverse of the exponent $d/(d+2)$ is refered as {\em Magnification 
Factor}. Note that in any case, when the data dimension is large, this 
exponent is near 1 (it value is 1/3 when $d=1$). Note also that this power 
has no effect when the density $f$ is uniform. But in fact the optimal 
quantizer is another thing, with another definition. 

\medskip

Namely the optimal quantizer $m^{\star}_n$ (formed with the code vectors 
$m^{\star}_{1,n}, m^{\star}_{2,n}, \ldots, m^{\star}_{n,n}$), minimizes the 
distortion $V_n(m)$, and is got after convergence of the 0-neighbor 
algorithm (if we could ensure the convergence to a global minimum, that is 
true only in the one-dimensional case). So if we set 
$$A_n(f,m^{\star}_n) = n^{\beta} V_n(m^{\star}_n) = 
      n^{\beta} \sum _{i=1}^{n} \int_{C_i} \|m^{\star}_{i,n}-x\|^2 
					f(x) dx $$

actually we have, 
$$A(f)=\lim_ {n \longrightarrow \infty} A_n(f,m^{\star}_n) < A(f,g^{\star}) $$

and the limit of the discrete distribution of $m^{\star}_n$ is not equal to 
$g^{\star}$. {\em So there is no magnification factor, for the 0-neighbor 
algorithm as claimed in many papers. It can be an approximation, but no more.}

The problem comes from the confusion between two distinct notions: 
random quantizer and optimal quantizer. And in fact, the good property 
is the convergence of the weighted distribution function (\ref{VQD}).


As to the SOM algorithm in the one-dimensional case, with a neighborhood 
function not reduced to the 0-neighbor case,  one can find in \cite{RITT5} 
or \cite{DERS} some result about a possible limit of the discrete 
distribution when the number of units goes to $\infty$. But actually, the 
authors use the Zador's result which is not appropriate as we just see.

\section{The multidimensional continuous setting}

In this section, we consider a general neighborhood function and the SOM 
algorithm is defined as in Section 2.

\subsection{Self-organization}

When the dimension $d$ is greater than 1, little is known on the classical 
Kohonen algorithm.  
The main reason seems to be the fact that it is difficult to define what can 
be an {\em organized state} and that no {\em absorbing} sets have been found. 
The configurations whose coordinates are monotoneous are not stable, 
contrary to the intuition. For each configuration set which have been claimed 
to be left stable by the Kohonen algorithm, it has been proved later that it 
was possible to go out with a positive probability. See for example 
\cite{COTT1}. Most people think that the Kohonen algorithm in dimension 
greater than 1 could correspond to an irreducible Markov chain, that is a 
chain for which there exists always a path with positive probability to go from 
anywhere to everywhere. That property  imply that there is no absorbing 
set at all. 

Actually, as soon as $d \geq 2$, for a constant parameter $\varepsilon$, 
the 0-neighbor algorithm is an Doeblin recurrent irreducible chain 
(see \cite{BOUT3}), that cannot have any absorbing class.

Recently, two apparently contradictory results were established, that can be 
collected together as follows.

\begin{monthe}

($d = 2$ and $\varepsilon$ is a constant) 
\noindent Let us consider a $n \times n$ units square network and the set $F^{++}$ 
of states whose both coordinates are separately increasing as function of 
their indices, i.e.
$$F^{++} = \left\{\forall i_1 \leq n, 
m_{i_1,1}^2 < m_{i_1,2}^2< \ldots < m_{i_1,n}^2, 
\forall i_2 \leq n,
m_{1,i_2}^1 < m_{2,i_2}^1 < \ldots < m_{n,i_2}^1 \right\}$$

\noindent (i) If $\mu$ has a density on $\Omega$, and if the neighborhood 
function $\Lambda$ is everywhere positive and decreases with the distance, 
the hitting time of $F^{++}$ is finite with positive probability (i.e. $>0$, 
but possibly less than 1). See Flanagan (\cite{FLAN}, \cite{FLAN1}).

\noindent (ii) In the 8-neighbor setting, the exit time from $F^{++}$ is 
finite with positive probability. See Fort and Pag\`es in  (\cite{FORT6}).
\label{orga2}
\end{monthe}

\medskip

This means that (with a constant, even very small, parameter $\varepsilon$), 
the organization is temporarily reached and that even if we guess that it is 
almost stable, dis-organization may occur with positive probability.

More generally, the question is how to define an organized state. Many 
authors have proposed definitions and measures of the self-organization, 
\cite{ZREH}, \cite{DEMA}, \cite {VILL}, \cite{GOOD}, \cite{VILL2}, 
\cite{HERR}. But none such ``organized'' sets have a chance to be absorbing. 

\medskip

In \cite{FORT6}, the authors propose to consider that {\em a map is 
organized if and only if the Vorono\"{\i} classes of the closest neighboring 
units are contacting}. They also precisely define the nature of the 
organization (strong or weak). 

They propose the following definitions :

\begin{madef} {\bf Strong organization}\\ 
There is strong organization if there exists a set of organized states 
${\cal S}$ such that

\noindent (i) ${\cal S}$ is an absorbing class of the Markov chain $m(t)$,

\noindent (ii) The entering time in ${\cal S}$ is almost surely finite, 
starting from any random weight vectors (see \cite{BOUT2}).
\end{madef}

\begin{madef} {\bf Weak organization}\\
There is weak organization if there exists a set of organized states 
${\cal S}$ such that all the possible attracting equilibrium points of the 
ODE defined in \ref{ODE} belong to the set ${\cal S}$.
\end{madef}

The authors prove that there is no strong organization at least in two 
seminal cases: the input space is $[0,1]^2$, the network is 
one-dimensional with two neighbors or two-dimensional with eight neighbors. 
The existence of weak organization should be investigated as well, but 
until now no exact result is available even if the simulations show 
a stable organized limit behavior of the SOM algorithm.

\subsection{Convergence}

In \cite{FORT5}, (see also \cite{FORT4}) the gradient of $h$ is 
computed in the $d$-dimensional setting (when it exists). In \cite{RITT2}, 
the convergence and the nature of the limit state is studied, assuming that 
the organization has occured, although there is no mathematical proof of 
the convergence. 

\medskip

Another interesting result received a mathematical proof thanks to the 
computation of the gradient of $h$: it is the dimension selection effect 
discovered by Ritter and Schulten (see \cite{RITT2}). The mathematical 
result is (see \cite{FORT5}:

\begin{monthe}

Assume that $m_1^{\star}$ is a stable equilibrium point of a general 
$d_1$-dimensional Kohonen algorithm, with $n_1$ units, stimuli distribution 
$\mu_1$ and some neighborhood function $\Lambda$. Let $\mu_2$ be a 
$d_2$-dimensional distribution with mean $m_2^\star$ and covariance matrix 
$\Sigma_2$. Consider the $d_1+d_2$ Kohonen algorithm with the same units and 
the same neighborhood function. The stimuli distribution is now 
$\mu_1 \bigotimes \mu_2$.

Then there exists some $\eta  > 0$, such that if $\| \Sigma_2 \| < \eta$,
the state $m_1^\star$ in the subspace $m_2 = m_2^\star$ is still a stable 
equilibrium point  for the $d_1 + d_2$ algorithm.
\label{variance}
\end{monthe}

It means that if the stimuli distribution is close to a $d_1$-dimensional 
distribution in the $d_1 + d_2$ space, the algorithm can find a $d_1$-space 
stable equilibrium  point. That is the {\em dimension selection effect}.

\medskip

From the computation of the gradient $\nabla h$, some partial results on 
the stability of grid equilibriums can also be proved:

\medskip

Let us consider $I = I_1 \times I_2 \times \ldots \times I_d$ a 
$d$-dimensional array, with $I_l=\{1,2, \ldots, n_l\}$, for 
$1 \leq l \leq d$. Let us assume that the neighborhood function is a product
function (for example 8 neighbors for $d=2$) and that the input 
distributions in each coordinate are independent, that is 
$\mu = \mu_1 \bigotimes \ldots \bigotimes \mu_d$. At last 
suppose that the support of each ${\mu}_l$ is [0,1]. 

Let us call {\em grid states} the states 
$m^{\star} = (m_{i_ll}^{\star}, 1 \leq i_l \leq n_l, 1 \leq l \leq d)$, 
such that for every 
$1 \leq l \leq d$, $(m_{i_ll}^{\star}, 1 \leq i_l \leq n_l)$ 
is an equilibrium for the one-dimensional algorithm. Then the
following results hold \cite{FORT5} :

\begin{monthe}

\noindent (i) The grid states are equilibrium points of the ODE (\ref{ODE})
in the $d$-dimensional case.

\noindent (ii) For $d=2$, if $\mu_1$ and $\mu_2$ have strictly positive 
densities $f_1$ and $f_2$ on $[0, 1]$, if the neighborhood functions are 
strictly decreasing, the  grid equilibrium points are not stable as soon as 
$n_1$ is large enough and the ratio $\frac{n_1}{n_2}$ is large (or small) 
enough (i.e. when $n_1 \longrightarrow +\infty$ and $\frac{n_1}{n_2} 
\longrightarrow +\infty$ or $0$, see \cite{FORT5}, Section 4.3).


\noindent (iii) For $d=2$, if $\mu_1$ and $\mu_2$ have strictly positive 
densities $f_1$ and $f_2$ on $[0, 1]$, if the neighborhood functions are 
degenerated (0 neighbor case), $m^\star$ is stable if $n_1$ and $n_2$ are 
less or equal to 2, is not stable in any other case (may be excepted when 
$n_1 = n_2 = 3$).
\label{grids}
\end{monthe}

\medskip

The (ii) gives a negative property for the non square grid which can be 
related with  this one: the product of one-dimensional quantizers is not the 
correct vectorial quantization. But also notice that we have no result about 
the simplest case: the square grid equilibrium in the uniformly distributed 
case. Everybody can observe by simulation that this square grid is stable 
(and probably the unique stable ``organized'' state). Nevertheless, even if 
we can numerically verify that it is stable, using the gradient formula it 
is not mathematically proved even with two neighbors in each dimension! 

Moreover, if the distribution $\mu_1$ and $\mu_2$ are not uniform, generally 
the square grids are not stables, as it can be seen experimentally.

\section{The discrete case}

In this case, there is a finite number $N$ of inputs and 
$\Omega = \{ x_1, x_2, \, \ldots , \, x_N \}$. The input 
distribution is uniform on $\Omega$ that is 
$\mu(dx) = \frac{1}{N}\sum_{l=1}^{N}\delta_{x_l}$. It is the 
setting of many practical applications, like Classification or Data Analysis.

\subsection{The results}

The main result (\cite{KOHO5}, \cite{RITT4}) is that for not time-dependent 
general neighborhood, the algorithm locally derives from the potential

\begin{eqnarray*}
V_n(m)&=& \frac{1}{2N} \sum_{i=1}^{n} \sum_{x_l \in C_i(m)} 
(\sum_{j=1}^{n} \Lambda(i-j) \|  m_j - x_l \| ^2) \\
      &=&\frac{1}{2}\sum_{i=1}^{n}\int_{C_i(m)}\sum_{j=1}^{n}  
\Lambda(i-j) \| m_j-x \|^2)  \mu(dx)\\
      &=&\frac{1}{2} \sum_{i,j=1}^n \Lambda(i-j) 
   		\int_{C_i(m)}\|m_j-x\|^2 \mu(dx).
\end{eqnarray*}

\medskip
When $\Lambda(i,j)=1$ if $i$ and $j$ are neighbors, and if ${\cal V}(j)$ 
denotes the neighborhood of unit $i$ in $I$, $V_n(m)$ also reads
$$V_n(m) = \frac{1}{2} \sum_{j=1}^n 
     \int_{\cup_{i \in {\cal V}(j)}C_i(m)} \|m_j - x\|^2 \mu(dx).$$
\medskip

$V_n(m)$ is an {\em intra-class variance extended to the neighbor classes} 
which is a generalization of the distortion defined in Section 4 for the 
0-neighbor setting. But this potential does have many singularities and its 
complete analysis is not achieved, even if the discrete algorithm can be 
viewed as a stochastic gradient descent procedure. In fact, there is a 
problem with the borders of the Vorono\"{\i} classes. The set of all these 
borders along the process ${\bf m}(t)$ trajectories has measure 0, but it is 
difficult to assume that the given points $x_l$ never belong to this set.

Actually the potential is the true measure of the self-organization. It 
measures both clustering quality and proximity between classes. Its study 
should provide some light on the Kohonen algorithm even in the continuous 
case.

When the stimuli distribution is continuous, we know that the algorithm is 
not a gradient descent \cite{ERWI2}. However the algorithm can be seen then 
as an approximation of the stochastic gradient algorithm derived from the 
function $V_n(m)$. Namely, the gradient of $V_n(m)$ has a non singular 
part which corresponds to the Kohonen algorithm and a singular one which 
prevents the algorithm to be a gradient descent.  

This remark is the base of many applications of the SOM algorithm as well 
in combinatorial optimization, data analysis, classification, 
analysis of the relations between qualitative classifying variables.

\subsection{The applications}

For example, in \cite{FORT1}, Fort uses the SOM algorithm with a close 
one-dimensional string, in a two dimensional space where are located $M$ 
cities. He gets very quickly a very good sub-optimal solution. See also the 
paper \cite{ANGE}.

The applications in data analysis and classification are more classical. The 
principle is very simple: after convergence, the SOM algorithm provides a 
two(or one)-dimensional organized classification which permit a low 
dimensional representation of the data. See in \cite{KOHO6} an impressive 
list of examples.

In \cite{COTT4} and \cite{COTT7}, an application to forecasting is presented 
from a previous classification by a SOM algorithm.

\subsection{Analysis of qualitative variables}

Let us define here two original algorithms to analyse the relations between 
qualitative variables. The first one is defined only for two qualitative 
variables. It is called KORRESP and is analogous to the simple classical 
Correspondence Analysis. The second one is devoted to the analysis of 
any finite number of qualitative variables. It is called KACM and is similar 
to the Multiple Correspondence Analysis. See \cite{COTT2}, \cite{COTT3}, 
\cite{COTT8} for some applications.

For both algorithms, we consider a sample of individuals and a number $K$ of 
questions. Each question $k, k=1,2,\ldots,K$ has $m_k$ possible answers 
(or modalities). Each individual answers each question by choosing one and 
only one modality. If $M=\sum_{1 \leq k \leq}m_k$ is the total number of 
modalities, each individual is represented by a row $M$-vector with 
values in ${0, 1}$. There is only one 1 between the 1st component and the 
$m_1$-th one, only one 1 between the $m_1+1$-th component and the 
$m_1+m_2$-th one and so on. 

In the general case where $M > 2$, the data are summarized into a Burt Table 
which is a cross tabulation table. It is a $M \times M$ symmetric matrix and 
is composed of $K \times K$ blocks, such that the $(k,l)$-block $B_{kl}$ 
(for $k \neq l$) is the $(m_k \times m_l)$ contingency table which crosses 
the question $k$ and the question $l$. The block $B_{kk}$ is a diagonal 
matrix, whose diagonal entries are the numbers of individuals who have 
respectively chosen the modalities $1, 2, \ldots,m_k$ for question $k$. In 
the following, the Burt Table is denoted by $B$.

In the case $M=2$, we only need the contingency table $T$ which crosses 
the two variables. In that case, we set $p$ (resp. $q$) for $m_1$ (resp. 
$m_2$). \\

\medskip

{\bf The KORRESP algorithm}\\

In the contingency table $T$, the first qualitative variable has $p$ levels 
and corresponds with the rows. The second one has $q$ levels and corresponds 
with the columns. The entry $n_{ij}$ is the number of individuals categorized 
by the row $i$ and the column $j$. From the contingency table, the matrix of 
relative frequencies ($f_{ij}=n_{ij}/(\sum_{ij}n_{ij})$) is computed.

Then the rows and the columns are normalized in order to have a sum equal 
to 1. The row profile $r(i), 1 \leq i \leq p$ is the discrete probability 
distribution of the second variable given that the first variable has 
modality $i$ and the column profile $c(j), 1 \leq j \leq q$ is the discrete 
probability distribution of the first variable given that the second 
variable has modality $j$. The classical Correspondence Analysis is a 
simultaneous weighted Principal Component Analysis on the row profiles and 
on the column profiles. The distance is chosen to be the $\chi^2$ distance. 
In the simultaneous representation, related modalities are projected into 
neighboring points. 

To define the algorithm KORRESP, we build a new data matrix ${\cal D}$ :
to each row profile $r(i)$, we associate the column profile $c(j(i))$ which 
maximizes the probability of $j$ given $i$, and conversely, we associate to 
each column profile $c(j)$ the row profile $r(i(j))$ the most probable given 
$j$. The data matrix ${\cal D}$ is the $((p+q) \times (q+p))$-matrix whose 
first $p$ rows are the vectors $(r(i),c(j(i)))$ and last $q$ rows are the 
vectors $(r(i(j)),c(j))$. The SOM algorithm is processed on the rows of this 
data matrix ${\cal D}$. Note that we use the $\chi^2$ distance to look for 
the winning unit and that we alternatively pick at random the inputs among 
the $p$ first rows and the $q$ last ones. After convergence, each modality 
of both variables is classified into a Vorono\"{\i} class. Related modalities 
are classified into the same class or into neighboring classes. This method 
give a very quick, efficient way to analyse the relations between two 
qualitative variables. See \cite{COTT2} and \cite{COTT9} for real-world 
applications.\\

\medskip

{\bf The KACM Algorithm}\\

When there are more than two qualitative variables, the above method does 
not work any more. In that case, the data matrix is just the Burt Table $B$. 
The rows are normalized, in order to have a sum equal to 1. At each step, 
we pick a normalized row at random according to the frequency of the 
corresponding modality. We define the winning unit according to the $\chi^2$ 
distance and update the weights vectors as usual. After convergence, we get 
an organized classification of all the modalities, where related modalities 
belong to the same class or to neighboring classes. In that case also, the 
KACM method provides a very interesting alternative to classical Multiple 
Correspondence Analysis.

The main advantages of both KORRESP and KACM methods are their rapidity and 
their small computing time. While the classical methods have to use several 
representations with decreasing information in each, ours provide only one 
map, that is rough but unique and permit a rapid and complete interpretation.
See \cite{COTT3} and \cite{COTT8} for the details and financial applications.

\section{Conclusion}

So far, the theoretical study in the one-dimensional case is nearly complete. 
It remains to find the convenient decreasing rate to ensure the ordering. 
For the multidimensional setting, the problem is difficult. It seems that 
the Markov chain is irreducible and that further results could come from 
the careful study of the Ordinary Differential Equation (ODE) and from the 
powerful existing results about the cooperative dynamical systems. 

On the other hand, the applications are more and more numerous, especially 
in data analysis, where the representation capability of the organized data 
is very valuable. The related methods make up a large and useful set of 
methods which can be substituted to the classical ones. To increase their 
use in the statistical community, it would be necessary to continue the 
theoretical study, in order to provide quality criteria and performance 
indices with the same rigour as for the classical methods.

\subsection*{Acknowledgements}
We would like to thank the anonymous rewiewers for their helpful comments.
 
\bibliographystyle{abbrv}

\end{document}